\documentclass{ifacconf}

\usepackage{graphicx}      
\usepackage{natbib}        

\usepackage{amsmath,amsfonts}
\usepackage{url}
\usepackage{multirow}
\def\R{{\mathbb R}}
\def\N{{\mathbb N}}

\begin{document}
\begin{frontmatter}

\title{Tractable semidefinite bounds of positive maximal singular values \thanksref{footnoteinfo}} 

\thanks[footnoteinfo]{This work was supported by the Tremplin ERC Stg Grant ANR-18-ERC2-0004-01 (T-COPS project), the FMJH Program PGMO (EPICS project), as well as the PEPS2 Program (FastOPF project) funded by AMIES and RTE.
This work has benefited from the European Union's Horizon 2020 research and innovation program under the Marie Sklodowska-Curie Actions, grant agreement 813211 (POEMA) as well as from the AI Interdisciplinary Institute ANITI funding, through the French ``Investing for the Future PIA3'' program under the Grant agreement n$^{\circ}$ANR-19-PI3A-0004.}

\author[First]{Victor Magron} 
\author[Second]{Ngoc Hoang Anh Mai} 
\author[Third]{Yoshio Ebihara}
\author[Fourth]{Hayato Waki}

\address[First]{LAAS-CNRS, Université de Toulouse, CNRS, IMT, Toulouse, France (e-mail: vmagron@laas.fr).}
\address[Second]{LAAS-CNRS, Université de Toulouse, CNRS,  Toulouse, France (e-mail: nhmai@laas.fr)}
\address[Third]{Graduate School of Information Science and Electrical Engineering, Kyushu University, 744 Motooka, Nishi-ku, Fukuoka819-0395, Japan (e-mail: ebihara@ees.kyushu-u.ac.jp)}
\address[Fourth]{Institute of Mathematics for Industry, Kyushu University, 744 Motooka, Nishi-ku, Fukuoka819-0395, Japan (e-mail: waki@imi.kyushu-u.ac.jp)}       

\begin{abstract}                
We focus on computing certified upper bounds for the positive maximal singular value (PMSV) of a given matrix.   
The PMSV problem boils down to maximizing a quadratic polynomial on the intersection of the unit sphere and the nonnegative orthant. 
We provide a hierarchy of tractable semidefinite relaxations to approximate the value of the latter polynomial optimization problem as closely as desired.
This hierarchy is based on an extension of P\'olya's representation theorem. 
Doing so, positive polynomials can be decomposed as weighted sums of squares of $s$-nomials, where $s$ can be a priori fixed ($s=1$ corresponds to monomials, $s=2$ corresponds to binomials, etc.).
This in turn allows us to control the size of the resulting semidefinite relaxations.
\end{abstract}


\end{frontmatter}

\section{Introduction}
\label{sec:intro}
Let $\N$, $\R$, $\R_+$ be the sets of nonnegative integers,  real numbers and nonnegative real numbers, respectively.
In the present paper, we focus on the problem of computing the positive maximal singular value (PMSV) of a given real matrix $M \in \R^{n\times n}$, denoted by $\sigma_+(M)$.  
Providing PMSV bounds is crucial for certain induced norm analysis of discrete-time linear time-invariant systems, where the input signals are restricted to be nonnegative; see, e.g., \cite{ebihara2021l2}.
The squared value $\sigma_+(M)^2$ can be obtained by solving the following polynomial maximization problem:
\begin{equation}
\label{eq:PMSV}
\begin{array}{rl}
     \sigma_+(M)^2 & =  \sup\limits_{x\in\R^n_+}\{x^\top (M^\top M) x\,:\,\|x\|_2^2=1\} \,, \\
& = \sup\limits_{x\in\R^n_+}\{x^\top (M^\top M) x\,:\,\|x\|_2^2 \leq 1\}     \,.
\end{array}
\end{equation}

For $m \in \N \backslash \{0\}$, let $[m] := \{1,\dots,m\}$.
Note that \eqref{eq:PMSV} is a particular instance of a polynomial optimization problem (POP) on the nonnegative orthant:
\begin{equation}\label{eq:constrained.problem.poly.even}
\begin{array}{l}
f^\star:=\sup\limits_{x \in S} f( x)\,,
\end{array}
\end{equation}
where $f$ is a polynomial and $S$ is a basic closed semialgebraic set, i.e., the intersection of finitely many polynomial inequalities as follows:
\begin{equation}\label{eq:semial.set.def.2.even}
    S:=\{x\in\R^n\,:\,x_i\ge 0\,,\,i\in[n]\,,\,g_j(x)\ge 0\,,\,j\in[m]\}\,,
\end{equation}
for some $g_j\in\R[x]$, $j\in[m]$ with $g_m:=1$.
Letting $x^2:=(x_1^2,\dots,x_n^2)$ and $\check p(x):=p(x^2)$  whenever $p\in \R[x]$,  it follows that problem \eqref{eq:constrained.problem.poly.even} is equivalent to solving
\begin{equation}\label{eq:new.pop}
    f^\star=\sup_{x\in \check S}\check f\,,
\end{equation}
with
\begin{equation}\label{eq:new.semi.set}
    \check S:=\{x\in\R^n\,:\,\check g_j(x)\ge 0\,,\,j\in[m]\}\,.
\end{equation}
In the case of PMSV, one has $m=2$, $\check g_1 = 1 - \sum_{i\in[n]} x_i^4$, $\check g_2 = 1$ and $\check f = (x^2)^\top M^\top M x^2$.

In \cite{dickinson2015extension}, the authors state a specific constrained version of P\'olya's Positivstellensatz. 
Explicitly, if $f,g_1,\dots,g_m$ are homogeneous polynomials, $S$ is defined as in \eqref{eq:semial.set.def.2.even}, and $f$ is positive on $S\backslash\{0\}$, then  
    \begin{equation}\label{eq:PuVa.perturb}
    \begin{array}{l}
         (\sum_{j\in[n]}x_j)^kf=\sum_{j\in[m]}\sigma_jg_j\,,
    \end{array}
    \end{equation}
    for some $k \in \N$ and homogeneous polynomials $\sigma_j$ 
    with positive coefficients. 
They also construct a hierarchy of linear relaxations associated with \eqref{eq:PuVa.perturb} converging from above to $f^\star$.

\paragraph*{\textbf{Contributions}.}
Here, we extend this result to the case where the input polynomials $f$, $g_j$ are all even, i.e., invariants under any variable sign flip, to get a representation similar to \eqref{eq:PuVa.perturb}, where each multiplier $\sigma_j$ is a sum of $s$-nomial squares. 
The case $s=1$ corresponds to monomials, yielding a hierarchy of linear programming relaxations similar to  \cite{dickinson2015extension}, the case $s=2$ corresponds to binomials, yielding a hierarchy of second-order conic programming relaxations. 
Fixing $s>2$ in advance leads to a hierarchy of SDP relaxations  with controlled sizes and allows us to overcome the potential ill-conditioning of the linear relaxations. 
Our optimization framework can be applied in particular to the PMSV problem and we illustrate the related numerical performance at the end of the paper. 

\paragraph*{\textbf{Related works}.}
A well-known method to approximate $f^\star$ is to consider the hierarchy of SDP relaxations by \cite{lasserre2001global} based on the representation by \cite{putinar1993positive}. Namely if $f$ is positive on $S$ and $S$ involves the polynomial $g_1 = L - \|x\|_2^2$ for some $L > 0$, then $f = \sum_{j\in[m]} \sigma_j g_j$ for some SOS polynomials $\sigma_j$. 
The resuling hierarchy of SDP relaxations to approximate $f^\star$ from above is:
\begin{equation}\label{eq:sdp.Pu}
    \begin{array}{l}
         \inf \left\lbrace \lambda : \lambda -  f = \sum_{j\in[m]}\sigma_jg_j\,, \deg(\sigma_jg_j) \leq 2 k  \right\rbrace \,.
    \end{array}
    \end{equation}
Increasing $k$ imrpoves the accuracy of the resulting upper bounds but the size of the SDP becomes rapidly intractable. 
To overcome this computational burden, a remedy consists of exploiting sparsity of the input polynomials; see the work by \cite{waki2006sums} and \cite{chordaltssos, wang2020cs}.
Another research direction is to restrict each multiplier $\sigma_j$ to be a sum of $s$-nomials. 
When $s\in\{1,2\}$, we retrieve the framework by \cite{ahmadi2019dsos} based on scaled diagonally dominant sums of squares (SDSOS), generalized later on by \cite{gouveia2022sums} for arbitrary larger $s$.
However, the resulting hierarchy of linear or second-order conic programming relaxations does not produce a sequence of upper bounds converging to $f^\star$.
It turns out that convergence can be obtained by multiplying $f$ by the power of a given positive polynomial, e.g., $1+\|x\|_2^2$, or equivalently forcing the multipliers $\sigma_j$ to have  denominators. 
Such representations have been derived in \cite{polya1928positive,reznick1995uniform,putinar1999positive} and the latter led to the optimization framework  presented in \cite{mai2021positivity}. 

\section{Extending P\'olya's representation}
\label{sec:polya}
We now state our main result, which extends P\'olya's representation to even input polynomials. 
In addition, we provide a degree bound for the multipliers.
\begin{thm}\label{coro:compact.even}
Let $g_1,\dots,g_m$ be even polynomials such that $g_1 =L-\|x\|_2^2$ for some $L>0$ and $g_m = 1$.
Let $S$ be the semialgebraic set defined by
\begin{equation}
    S:=\{x\in\R^n\,:\,g_1(x)\ge 0\dots,g_m(x)\ge 0\}\,.
\end{equation} 
Let $f$ be an even polynomial and nonnegative on $S$ and $d_f: = \lfloor \deg(f)/2\rfloor +1$.
Then the following statements hold:
\begin{enumerate}
    \item For all $\varepsilon>0$, there exists $K_\varepsilon\in\N$ such that for all $k\ge K_\varepsilon$, there exist sums of monomial squares $\sigma_j$ with $\deg(\sigma_j g_j)\le 2(k+d_f)$, $j \in [m]$ and 
\begin{equation}\label{eq:reperese.compact}
   (1+\|x\|_2^2)^{k}(f+\varepsilon)=\sum_{j\in[m]} \sigma_j g_j \,.
\end{equation}
\item If $S$ has nonempty interior, there exist positive constants $\bar c$ and $c$ depending on $f$ and $S$ such that for all $\varepsilon>0$,  \eqref{eq:reperese.compact} holds with $K_{\epsilon} = \bar c\varepsilon^{-c}$.
\end{enumerate}
\end{thm}
Theorem \ref{coro:compact.even} can be proved in the same way as \cite[
Corollary 2]{mai2021complexity}.
It is important to note that the theorem still holds if we replace the first constraint polynomial $g_1 = L - \|x\|_2^2$  by $L - \sum_{i\in[n]} x_i^4$ and we do so for the PMSV problem.
If we remove the multiplier $(1+\|x\|_2^2)^k$ in \eqref{eq:reperese.compact}, Theorem \ref{coro:compact.even} is no longer true.
Indeed, let $n=1$, $f:=(x^2-\frac{3}{2})^2$ and assume that  $f=\sigma_1(1-x^2) + \sigma_2$ for some SOS of monomials $\sigma_1, \sigma_2$.
Note that $f$ is even and positive on $[-1,1]$.
We write $\sigma_i:=a_i+b_i x^2+x^4r_i(x)$ for some $a_i,b_i\in\R_+$ and $r_i\in\R[x]$.
It implies that
\begin{equation}
\begin{array}{rl}
     x^4-3x^2+\frac{9}{4}= &(a_1+b_1 x^2+x^4r_1(x))(1-x^2) \\
     & + (a_2+b_2 x^2+x^4 r_2(x))\,.
\end{array}
\end{equation}
Then we obtain the system of linear equations:
$\frac{9}{4}=a_1+a_2$ and $-3=b_2-a_1+b_1$.
Summing gives $-\frac{3}{4}=a_2+b_2+b_1$. 
However, $a_2+b_2+b_1\ge 0$ since $a_i,b_i\in\R_+$, yielding a contradiction.
Thus, Putinar's representation with sums of monomial squares does not exist for even input polynomials.
%
With the multiplier $(1+x^2)^2$, we obtain a P\'olya's representation:
\begin{equation}\label{eq:Pu-Va.decomp.ins}
    \begin{array}{rl}
         (1+x^2)^2f=  \bar \sigma_1(1-x^2) + \bar \sigma_2\,,
    \end{array}
\end{equation}
where  $\bar \sigma_1=x^4+\frac{15}{4}x^2+\frac{9}{4}$ and $\bar\sigma_2=x^8$ are SOS of monomials.

\section{Tractable semidefinite relaxations}
\label{sec:sdp}
Given $\alpha = (\alpha_1,\dots,\alpha_n) \in \N^n$, we write $x^{\alpha} = x_1^{\alpha} \cdots x_n^{\alpha_n}$.
An $s$-nomial square is a polynomial which can be written as $(\sum_{i\in[s]} a_i x^{\alpha(i)} )^2$, with $a_i \in \R$, $\alpha(i) \in \N^n$, $i \in [s]$.
Let us denote by $\Sigma_s$ the set of sums of $s$-nomial squares. 
The set of monomial squares corresponds to $\Sigma_1$.
For any $s \in \N \backslash \{0\}$, one has the obvious inclusion $\Sigma_1 \subset \Sigma_s$, thus the representation result \eqref{eq:reperese.compact} from Theorem \ref{coro:compact.even} holds with multipliers $\sigma_j \in \Sigma_s$, $j \in [m]$.
With $f$ and $S$ as in Theorem \ref{coro:compact.even}, we rely on 
this representation to derive a converging sequence of upper bounds for $f^\star$.  
Each upper bound is obtained by solving an SDP, indexed by  $s\in\N \backslash \{0\}$ and $k\in\N$: 
\begin{equation}\label{eq:primal.problem.0.even}
\begin{array}{rl}
   \rho _{k,s}:= \inf\limits_{\lambda} & \lambda\\
   \text{s.t.} & (1+\|x\|_2^2)^k (\lambda - f)=\sum_{j\in[m]} \sigma_j g_j \,, \\
   & \deg(\sigma_j g_j) \leq 2(k+d_f) \,, \  \sigma_j \in \Sigma_s  \,, j \in [m] \,.
\end{array}
\end{equation}
As a consequence of Theorem \ref{coro:compact.even}, we obtain the following convergence result.
\begin{cor}\label{theo:constr.theo.0.even}
Let $f$ and $S$ be as in Theorem \ref{coro:compact.even}. 
The following statements hold:
\begin{enumerate}
\item For all $k\in\N$ and for every $s\in\N \backslash \{0\}$,
    $\rho_{k,1} \ge \rho_{k,s} \ge f^\star$.
\item For every $s \in \N \backslash \{0\}$, the sequence $(\rho_{k,s})_{k\in\N}$ converges to $f^\star$.
\item If $S$ has nonempty interior, there exist positive constants $\bar c$ and $c$ depending on $f$ and $S$ such that for every $s\in\N \backslash \{0\}$ and for every $k\in\N$,
\begin{equation}
    0\le \rho_{k,s} - f^\star \le \left(\frac{k}{\bar c}\right)^{-\frac{1}{c}} \,.
\end{equation}
\end{enumerate}
\end{cor}

\section{Numerical experiments}
\label{sec:bench}
We demonstrate the accuracy and efficiency of our optimization framework, namely the SDP relaxations \eqref{eq:primal.problem.0.even} based on the extension of P\'olya's representation theorem.
Our framework is implemented in Julia 1.3.1 and available online via the link \url{https://github.com/maihoanganh/InterRelax}.
We compare the upper bounds and runtime with the ones of the standard SDP relaxations \eqref{eq:sdp.Pu}, based on Putinar's  Positivstellensatz.
These latter relaxations are modeled by TSSOS \cite{wang2021tssos}.  All SDP relaxations are solved by the solver Mosek 9.1. 
We use a desktop computer with an Intel(R) Core(TM) i7-8665U CPU @ 1.9GHz $\times$ 8 and 31.2 GB of RAM. 
Our benchmarks come from induced norm analysis of linear time-invariant systems, where the input signals are restricted to be nonnegative; we generate random matrices $M$ as in
\cite[(12)]{ebihara2021l2}, that is
\begin{equation}
    M:=\begin{bmatrix}
    D & 0 &0& \dots& 0\\
    CB& D &0&\dots& 0\\
    CAB& CB & D& \dots& 0\\
    \dots&\dots&\dots&\dots&\dots\\
    CA^{r-2}B&CA^{r-3}B&CA^{r-4}B & \dots &D
    \end{bmatrix}\,,
\end{equation}
where $A,B,C,D$ are square matrices of size $r$, fixed in advance.
In our numerical experiments, we choose $r \in \{4,5,6,7\}$ and every entry of $A,B,C,D$ is taken uniformly in $(-1,1)$.
\if{
In order to compute the PMSV $\sigma_+(M)$ of $M$, we solve the following POP on the nonnegative orthant:
\begin{equation}
\begin{array}{rl}
     \sigma_+(M)^2=\max\limits_{x\in\R^n_+}\{x^\top (M^\top M) x\,:\,\|x\|_2^2=1\}\,.
\end{array}
\end{equation}
}\fi
Note that the number of variables of the resulting POP is $n=r^2$. 
We associate  an ``Id'' to each related SDP relaxation, in which we compute a decomposition into sums of $s$-nomial squares.
For Putinar's representation, we take a small relaxation order $k \in \{1,2\}$ while for the extension P\'olya's representation we let $k=0$,  leading to decompositions without denominators. 
We indicate the data of each SDP program, namely ``nmat'' and  ``msize''  correspond to the number of matrix variables and their largest size, ``nscal'' and ``naff'' are the number of scalar variables and affine constraints. 
The symbol ``$-$'' means that the corresponding computation aborted due to lack of memory.
The most accurate upper bounds are emphasized in bold. 
The total running time (including modeling and solving time) is indicated in seconds. 

The numerical results are displayed in Table \ref{tab:PMSV}.
The SDP relaxations associated to our extension of P\'olya's representation provide more accurate upper bounds and they are computed more efficiently.

\begin{table}[!ht]
    \caption{\small Upper bounds of PMSV of $M$.}
    \label{tab:PMSV}
\tiny
\begin{center}
\begin{tabular}{|c|c|c|c|c|c|c|c|c|}
        \hline
\multirow{2}{*}{Id}     & \multirow{2}{*}{$n$}
&
\multicolumn{3}{c|}{Putinar}&      \multicolumn{4}{c|}{P\'olya}\\ 
\cline{3-9}
&
&
{$k$} & {upper bound} & {time} & {$k$} & {$s$}& {upper bound} & {time}\\
\hline
\hline
1 & \multirow{2}{*}{16} &  1 & 47.48  & 0.02& \multirow{2}{*}{0} & \multirow{2}{*}{17} &  \multirow{2}{*}{\bf 30.18} & \multirow{2}{*}{0.7} \\
2 &  &  2 & {\bf 30.18} & 16&    & &  & \\
\hline
\hline
3 & \multirow{2}{*}{25} &  1 & 168.44  & 0.04& \multirow{2}{*}{0} & \multirow{2}{*}{26} &  \multirow{2}{*}{\bf 91.28}  & \multirow{2}{*}{0.9} \\
4&  &  2 & {\bf 91.28} & 877 & &  &  & \\
\hline
\hline
5 & \multirow{2}{*}{36} &  1 & 4759.12  & 0.2& \multirow{2}{*}{0} & \multirow{2}{*}{37} &  \multirow{2}{*}{\bf 2462.03}  & \multirow{2}{*}{1} \\
6& &  2 & $-$ & $-$&  &  &  & \\
\hline
\hline
7 & \multirow{2}{*}{49} &  1 & 1777.53  & 0.5& \multirow{2}{*}{0} & \multirow{2}{*}{50} &  \multirow{2}{*}{\bf 970.20}  & \multirow{2}{*}{2} \\
8&  &  2 & $-$ & $-$& & &  & \\
\hline
\end{tabular}

\vspace{0.5cm}

\begin{tabular}{|c|c|c|c|c|c|c|c|c|c|c|c|c|c|c|c|c|}
        \hline
\multirow{2}{*}{Id}&
\multicolumn{4}{c|}{Putinar}&      \multicolumn{4}{c|}{P\'olya}\\ \cline{2-9}
&
nmat & msize & nscal & naff &nmat & msize & nscal & naff\\
\hline
\hline
1 & 1 & 17 & 38 & 153 & \multirow{2}{*}{1} & \multirow{2}{*}{17} & \multirow{2}{*}{138} & \multirow{2}{*}{153}\\
2 & 17 & 153 & 154 & 4845 & &  & &\\
\hline
\hline
3 & 1 & 26 & 27 & 351 & \multirow{2}{*}{1} & \multirow{2}{*}{26} & \multirow{2}{*}{327} &\multirow{2}{*}{351}\\
4 & 26 & 351 & 352 & 23751 & &  & &\\
\hline
\hline
5 & 1 & 37 & 38 & 703 & \multirow{2}{*}{1} & \multirow{2}{*}{37} & \multirow{2}{*}{668} &\multirow{2}{*}{703}\\
6 & 37 & 703 &704 & 91390 &  & &&\\
\hline
\hline
7 & 1 & 50 & 51 & 1275 & \multirow{2}{*}{1} & \multirow{2}{*}{50} & \multirow{2}{*}{1227} &\multirow{2}{*}{1275}\\
8 & 50 & 1275 &1276 & 292825 & & & &\\
\hline
\end{tabular}
\end{center}
\end{table}


\end{document}